\documentclass[english]{amsart}
\usepackage[all]{xy}
\usepackage[dvips]{graphics}
\usepackage[latin1]{inputenc}
\usepackage{color}
\usepackage{amssymb}
\usepackage{amsmath}

\newcommand{\PrpTitle}[1]{{\upshape \bf  (}#1{\upshape\bf )}}

\newcommand{\im}{{\mathbf i}}

\newcommand{\adj}{{\rm Adj}}
\newcommand{\cre}{{\rm Bir}}

\newcommand{\dec}{{\rm Dec}}
\newcommand{\ine}{{\rm Ine}}
\newcommand{\aut}{{\rm Aut}}
\newcommand{\Aut}{{\rm Aut}}

\newcommand{\bir}{{\rm Bir}}
\newcommand{\Bir}{{\rm Bir}}

\newcommand{\pic}{{\rm Pic}}

\newcommand{\pgl}{{\rm PGL}}

\newcommand{\proj}{{\mathbb P}}

\newcommand{\plan}{{\mathbb P}^2}
\newcommand{\Pn}{{\mathbb P}^2}

\newcommand{\complex}{{\mathbb C}}

\newcommand{\C}{{\mathbb C}}

\newcommand{\Z}{{\mathbb Z}}

\newcommand{\hts}{\vphantom{\big(}}

\newcommand{\wtilde}{\widetilde}
\newcommand{\tor}{\xymatrix{\ar@{-->}[r]&}}

\newtheorem{thm}{Theorem}[subsection]
\newtheorem{pro}[thm]{Proposition}
\newtheorem{ques}[thm]{Question}
\newtheorem{exa}[thm]{Example}
\newtheorem{exas}[thm]{Examples}
\newtheorem{cor}[thm]{Corollary}
\newtheorem{lem}[thm]{Lemma}
\newtheorem{defi}[thm]{Definition}

\theoremstyle{remark}
\newtheorem{rem}[thm]{Remark}

\title[Birational Transformations of pairs]{On Birational Transformations of Pairs in the Complex Plane}
\author{J\'er\'emy Blanc, Ivan Pan, Thierry Vust}

\address{J\'er\'emy Blanc,
Universit\'e de Grenoble I,
UFR de Math\'ematiques,
UMR 5582 du CNRS, 
Institut Fourier, BP 74, 
38402 Saint-Martin d'H\`eres, France}
\address{
Ivan Pan, 
Instituto de Matem\'atica, 
UFRGS, 
av.~Bento Gon\c calves 9500, 
91540-000 Porto Alegre, 
RS, Brasil}
\address{Thierry Vust, 
Universit\'e de Gen\`eve,  
Section de math\'ematiques, 
2-4 rue du Li\`evre, 
CP 64, 
1211 Gen\`eve 4, 
Switerzland}
\begin{document}
\maketitle
\begin{abstract}
This article deals with the study of the birational transformations of the projective complex plane which leave invariant an irreducible algebraic curve.
We try to describe the state of art and provide some new results on this subject.
\end{abstract}

\section{Introduction}
\subsection{The decomposition and inertia groups}
We study the birational transformations of the projective complex plane which leave invariant an irreducible algebraic curve.

We  denote by 
$\Bir(\Pn)$ the group of birational transformations of the complex projective plane $\Pn=\Pn(\C)$: this is the \emph{Cremona group} of $\plan$.
If $C\subset \plan$ is an irreducible curve and $\varphi\in\Bir(\Pn)$, we say that $\varphi$ \emph{preserves} $C$ (or \emph{leaves $C$ invariant}) if $\varphi$ restricts to a birational transformation of $C$. If this latter is the identity, we say that $\varphi$ \emph{fixes}~$C$.

Let $C\subset\plan$ be an irreducible plane curve. 
Following Gizatullin \cite{Giz}, we introduce the \emph{decomposition group} of $C$ in $\Bir(\Pn)$, here denoted by $\dec(\Pn,C)=\dec(C)$, as the group of Cremona transformations that preserve $C$. The action $\rho$ of $\dec(C)$ on $C$  induces a (not necessarily exact) complex

\begin{equation}\xymatrix{1\ar@{->}[r]&\ine(C)\ar@{->}[r]&\dec(C)\ar@{->}[r]^\rho&\bir(C)\ar@{->}[r]&1,}\label{CanComplP2}\end{equation}
where $\ine(C)=\ine(\Pn,C):=\ker(\rho)$ is the \emph{inertia group} of $C$ in $\Bir(\Pn)$, which is the group of Cremona transformations that fix $C$. 
\subsection{Birational geometry of pairs}
The above notions may in fact be generalised to pairs $(S,C)$, where $S$ is a surface and $C\subset S$ an irreducible curve. We say that a birational transformation $\varphi:S\dasharrow S'$ is a \emph{birational transformation of pairs} $\varphi:(S,C)\dasharrow (S',C')$ if it restricts to a birational transformation $\varphi_{|_C}:C\dasharrow C'$, and in this case we say that the two pairs are \emph{birationally equivalent}. The group of birational transformations  of a pair $(S,C)$ is denoted by $\dec(S,C)$ and induces as before a complex 
\begin{equation}\xymatrix{1\ar@{->}[r]&\ine(S,C)\ar@{->}[r]&\dec(S,C)\ar@{->}[r]^\rho&\bir(C)\ar@{->}[r]&1},\label{CanComplGen}\end{equation}
which is exactly the complex (\ref{CanComplP2}) if $S=\Pn$ (and in this case we remove the surface in the notation). We will say that (\ref{CanComplGen}) is the \emph{canonical complex} of the pair $(S,C)$. Note that $\Aut(S,C):=\Aut(S)\cap \dec(S,C)$ is the group of automorphisms of $S$ that leave invariant the curve $C$.

\subsection{Outline of the article}
The aim of this article is to give a survey about the pairs $(\plan,C)$ whose decomposition group is not trivial together with a description of their corresponding canonical complexes; we point out what is known to us about the subject and we give some new results. 

Sections~\ref{Sec:GenusAtLeast2}, \ref{Sec:Genus1} and~\ref{Sec:Genus0} deal respectively with the curves of genus $\geq 2$, $1$ and $0$. Sections~\ref{Sec:LinkFinite} and~\ref{Sec:Dynamic} express the link between the transformations that preserve or fix curves with respectively the classification of finite subgroups and the dynamic of the elements of $\Bir(\Pn)$.

\subsection{Conventions}
In the sequel, $g(C)$ will denote the geometric genus of an irreducible curve $C$. Recall also that a \emph{de Jonqui\`eres transformation} is a birational transformation of $\Pn$ that preserves a pencil of lines. Finally, all our surfaces will be assumed to be rational, smooth, projective and irreducible.
\section{Curves of genus at least equal to $2$}\label{Sec:GenusAtLeast2}
\subsection{The main tool: adjoint linear systems}
Let $C\subseteq \plan$ be an irreducible curve with $g(C)\geq 2$. To study the group $\dec(C)$, we follow an idea of Castelnuovo and Enriques which consists of considering the \emph{adjoint} linear system associated to $C$: we take an embedded resolution of the singularities of $C$, say $\sigma:Y\to \plan$, denote by $\tilde{C}\subset Y$ the strict transform of $C$ and consider the linear system $\sigma_*|K_Y+\wtilde{C}|$ and its fixed part $\Delta$. By definition, the adjoint system $\adj(C)$ is the linear system $\sigma_*|K_Y+\wtilde{C}|-\Delta$. By Riemann-Roch it has dimension $g(C)-1>0.$ The main result is:

\begin{pro}
\PrpTitle{\cite{Cas}, \cite[Prop. 2.5]{BPV}}

Let $C\subset \Pn$ be an irreducible curve with $g(C)\geq 2$.

If $\varphi\in\bir(\Pn)$ sends $(\Pn,C)$ on $(\Pn,D)$, then it sends $\adj(C)$ on $\adj(D)$. In particular, the group $\dec(C)=\dec(\Pn,C)$ stabilizes the linear system $\adj(C)$. 
\end{pro}
One can also define the adjoint of a linear system, by taking the adjoint of a general member of the system. Since this construction decreases the degree of the curves, it has to stop after a finite number of iterations when the curves have no adjoint, i.e. when the curves have genus $0$ or $1$. This yields the following:
\begin{pro}\label{Prp:InvariantPencil}
\PrpTitle{\cite[Prop. 2.12]{BPV}}

Let $C\subset \Pn$ be an irreducible curve, with $g(C)\geq 2$. 
There exists a linear system (resp.~a pencil) of elliptic or rational curves $\Lambda$ such that $\dec(C)$ (resp.~$\ine(C)$) stabilizes $\Lambda$.
\end{pro}

\subsection{The inertia group of curves of genus $\geq 2$}
G. Castelnuovo used the existence of the invariant pencils yielded by Proposition~\ref{Prp:InvariantPencil} to bound the order of elements of finite order of $\ine(C)$:

\begin{thm}\PrpTitle{\cite{Cas}, \cite[Chap.VIII, \S 2]{God}, \cite[Book IV, Chap. VII, \S 3]{Coo}}
Let $C\subset \Pn$ be an irreducible curve with $g(C)\geq 2$, and let $\varphi\in \ine(C)$, $\varphi\not=1$. Then, either $\varphi$ is conjugate to a de Jonqui\`eres transformation or $\varphi$ has order $2,3$ or $4$. 
\end{thm}

In \cite{BPV}, an examination of the two possible cases of pencils yielded by Proposition~\ref{Prp:InvariantPencil} leads to a precise description of all cases of pairs $(\Pn,C)$ having a non-trivial inertia group. This generalises Castelnuovo's Theorem. We describe these cases in Examples~\ref{Exa:GeiBerJon}, \ref{Exa:dP1Ordr3} and~\ref{Exa:deJAbGr} below and then state the classification (Theorem~\ref{THM:BPVClass}). 



\begin{exas}\label{Exa:GeiBerJon} 
\PrpTitle{\cite{Hud}, \cite{God}, \cite{Coo}, \cite{SR}, \cite{BayBea},
\cite{Fer}, \cite{BlaCR}}

a) Let $p_1,...,p_7$ be $7$ points in the plane in general position. The Geiser involution does the following: a general point $q$ in the plane defines a pencil of cubic curves passing through $q$ and the seven points $p_1,...,p_7$; this pencil has a ninth base-point, which is the image of $q$ by the Geiser involution. This involution fixes a non hyperelliptic curve of genus 3 that is a sextic with ordinary double points at $p_1,...,p_7$ and whose smooth model is a plane quartic {\upshape (}\cite{Hud} or \cite{God}{\upshape )}; the blow-up of the seven points conjugates the Geiser transformation to an automorphism of a del Pezzo surface of degree~$2$ {\upshape (}\cite{BayBea}{\upshape )}.

b) Let $p_1,...,p_8$ be $8$ points in the plane in general position. The Bertini involution does the following: a general point $q$ in the plane defines a linear system of sextic curves passing through $q$ and  being singular at each of the eight points $p_1,...,p_8$. This linear system has a tenth base-point, which is the image of $q$ by the Bertini involution. This involution fixes a non hyperelliptic curve of genus 4 that is a nonic with ordinary triple points at $p_1,...,p_8$ and whose smooth model lies on a quadratic cone; the blow-up of the eight points conjugates the Bertini transformation to an automorphism of a del Pezzo surface of degree~$1$ (same references as above).

c) Let $C\subset \Pn$ be a curve of degree $g+2$ with an ordinary $g$-fold point and which is smooth anywhere else. The de Jonqui\`eres  involution associated to $C$ is defined in the following way: the restriction of the transformation to a general line passing through the $g$-fold point of $C$ is the unique involution that preserves this line and fixes the two other points of intersection of $C$ with the line {\upshape (}\cite{Jon}, \cite{BayBea}{\upshape )}. 
\end{exas}

\begin{exa}\label{Exa:dP1Ordr3} 
\PrpTitle{\cite{Fer}, \cite{Dol}, \cite{BlaCR}}

Consider the smooth surface $S$ defined by the equation  
$w^2=z^3+F_6(x,y)$ in the weighted projective space $\proj(3,1,1,2)$, where $F_6$ is an homogeneous polynomial of degree 6 with $6$ simple roots: it is a particular type of  del Pezzo surface of degree 1. The restriction of the map 
$(w:x:y:z)\mapsto (w:x:y:\omega z)$, where $\omega\neq 1$ is a primitive cubic root of 1, defines an automorphism of $S$ of order 3 whose set of fixed points is the union of a point and an irreducible curve $\tilde{C}$ of genus 2. The curve being linearly equivalent to $-2K_S$, any birational morphism $S\rightarrow \plan$ sends this curve on a sextic with $8$ ordinary double points in general position.\label{exa-order3}
\end{exa}

\begin{exa}\label{Exa:deJAbGr} 
\PrpTitle{\cite{BPV}}

Let $h\in\C[x]$ be a polynomial of degree $2g+2$ without multiple roots. Consider the subgroup 
\[T_h:=\Big\{\left(\begin{array}{ll} a_1&ha_2\\a_2&a_1\end{array}\right):
a_i\in\complex(x), a_1^2-h a_2^2\neq 0\Big\}\]
of ${\rm GL}(2,\complex(x))$ and denote by $J_h$ its image in $\pgl(2,\complex(x))$.
To each $a\in J_h$,  we associate a rational map 
$F_a:\complex^2\tor\C^2$ defined by
\[(x,y)\mapsto \left(x,\frac{a_1y+ha_2}{a_2y+a_1}\right):\]
it is a de Jonqui\`eres transformation whose restriction to the hyperelliptic curve
$C$ of equation $(y^2=h(x))$ is the identity. When $a_1=0$, we obtain an involution $\sigma$, conjugate to Example~\ref{Exa:GeiBerJon}c.

Notice that $T_h$ is isomorphic to the multiplicative group $\C(C)^{*}$  of the field of rational functions $\C(C)$ on $C$ from which we deduce that $J_h$ is isomorphic to $\C(C)^*/\C(x)^*$ and that its torsion is generated by $\sigma$.
\end{exa}

 
\begin{thm}\PrpTitle{\cite[Theorem 1.5]{BPV}}
Let $C\subset \Pn$ be an irreducible curve of genus $g\geq 2$, and assume that $\ine(C)$ is non-trivial. 

Then, either $\ine(C)$ is a cyclic group of order $2$ or $3$ generated by one of the transformations from Examples~\ref{Exa:GeiBerJon}a), \ref{Exa:GeiBerJon}b), \ref{Exa:dP1Ordr3} or it is equal the group $J_h$ of Example~\ref{Exa:deJAbGr}, where $(y^2=h(x))$ is the affine equation of $C$.

In particular, $\ine(C)$ is Abelian and if $\ine(C)$ is infinite, then $C$ is hyperelliptic and $\ine(C)$ is a de Jonqui\`eres group, whose torsion is generated by a de Jonqui\`eres involution.
\label{THM:BPVClass}
\end{thm}
Remark that Theorem~\ref{THM:BPVClass} implies in particular that the elements of order $4$ suggested by Castelnuovo's theorem do not exist.  It implies also the following result, proved in \cite[Theorem 1.1]{PanCR}.\begin{cor}
Let $S$ be a projective smooth rational surface and let $C\subset S$ be an irreducible curve with $g(C)>1$. Then, the group of elements of $\Aut(S)$ that fix $C$ has order $1$, $2$ or $3$. 
\end{cor}
\begin{proof}
Let us write $G=\ine(S,C)\cap \Aut(S)$. According to Theorem~\ref{THM:BPVClass}, we may assume that $G$ preserve a rational fibration $p:S\dasharrow \mathbb{P}^1$, and it suffices to show that no element of $G$ is of infinite order. Suppose for contrary that some $\varphi \in G$ is of infinite order. After some blow-up we may suppose that $p$ is a morphism (since $\varphi$ acts on the base-point of  the fibration). Then, we replace $\varphi$ by some power, and assume that $\varphi$ preserves any component of any singular fibre of $p$. This implies that $\varphi$ is conjugate to an automorphism of a Hirzebruch surface, which is not possible since it fixes a curve of positive genus.
\end{proof}

\subsection{The decomposition group of curves of genus $\geq 2$}
Applying the classification of the non-trivial inertia groups of curves of genus at least $2$, we deduce the following:

\begin{thm}
Let $C\subset \Pn$ be an irreducible curve of genus $g\geq 2$, and assume that $\ine(C)$ is non-trivial. Then, the canonical complex of $(\Pn,C)$ is an exact sequence.
\label{THM:SplitExactSequenceWhenIne1}
\end{thm}

\begin{proof}
Theorem~\ref{THM:BPVClass} reduces the possibilities for the pair $(\Pn,C)$.

The exactness of the canonical complex in the case where $\ine(C)$ is generated by the Geiser or Bertini involution is classical. For a proof (see \cite[Thm. 1.8]{Pan}), we consider the decomposition group as a subgroup of automorphisms of a del Pezzo surface $S$ of degree $2$ or $1$ and denote by $\sigma$ the Geiser or the Bertini involution; then each automorphism of the curve $\Gamma$ fixed by $\sigma$ extends to an automorphism of $S$ because $\Gamma$ is canonical in $S/<\sigma>$. 

In the de Jonqui\`eres case (Examples~\ref{Exa:GeiBerJon}c) and~\ref{Exa:deJAbGr}), denote by $C$ the curve of degree $g+2$, by $p$ its $g$-fold point and by $\Lambda$ the pencil of lines passing by $p$. Let $j$ be a birational map of $C$. Since $C$ is hyperelliptic $j$ stabilizes the trace of $\Lambda$ on $C$. Let $x\in\plan$ a general point. We extend $j$ to an element $\chi$ of $\dec(C)$: indeed, take the line $L\in\Lambda$ passing by $p$ and $x$ and set $L\cap C=\{p,q_x,r_x\}$; we define $\chi(x)$ by the relation
\[(p,x,q_x,r_x)=(p,\chi(x),j(q_x),j(r_x)),\] 
where $(a,b,c,d)$ denotes the cross ratio of $a,b,c,d$.

In the last case we consider the pair $(S,\widetilde{C})$ as in Example~\ref{exa-order3} and we observe that the restriction homomorphism $\aut(S,\widetilde{C})\to\aut(\widetilde{C})$ is surjective: indeed, an automorphism of $\widetilde{C}$ extends to an automorphism of $\proj(3,1,1)$ which  lifts to an automorphism of $S$.
\end{proof}

\begin{rem}
It may be observed that the exact sequence described above is split in the de Jonqui\`eres and Geiser cases, and in the case of Example~\ref{Exa:dP1Ordr3}. However it does not split in the Bertini case (there are square roots of the Bertini involution, see \cite[Table 1]{BlaCR}). 
\end{rem}

Theorem~\ref{THM:SplitExactSequenceWhenIne1} achieves the classification of pairs $(\plan, C)$ and canonical complexes such that $\ine(C)\neq 1$. 

What happens in the case where the group $\ine(C)$ is trivial? 
Firstly, we can state the following obvious result:

\begin{lem}
Let $C\subset \Pn$ be an irreducible curve of genus $g\geq 2$, and assume that $\ine(C)=1$. Then, $\dec(C)$ is isomorphic to a subgroup of $\cre(C)$, and is a finite group. In particular, when $C$ is generic,  $\ine(C)=\dec(C)=1$.
\end{lem}
\begin{proof}
It suffices to observe that $\cre(C)$ is isomorphic to the automorphisms group of the normalization of $C$ and recall that this group is finite when $g(C)\geq 2$ and is trivial if $C$ is generic.
\end{proof}
The canonical complex is therefore trivially exact for a general curve. However, there exist examples where the map $\dec(C)\rightarrow \cre(C)$ is not surjective, see Sections~\ref{SubSecGenus3}, \ref{SubSec:Genus1NotCubic}, \ref{SubSecCoble} and~\ref{SubSecRationalofHalphenType}. These rely on the following result:

\begin{thm}\label{Thm:BasePoints}
Let $C\subset \Pn$ be an irreducible curve of degree $n$. For each point $p$ that belongs to $C$ as a proper or infinitely near point, we denote by $m_p$ the multiplicity of $C$ (or its strict transform) at $p$. 

Let $\varphi\in \Bir(\Pn)$.
Suppose that $3m_p\leq n$ for each point $p$ and that $\varphi$ sends $C$ on a curve $D$ of degree $\leq n$.

For every base-point $q$ of $\varphi$ the point $q$ belongs to $C$ as a proper or infinitely near point, and $3m_q=n$. Moreover, the degree of $D$ is $n$.
\end{thm}
\begin{proof}We may assume that $\varphi$ is not an automorphism of $\Pn$. 
Let $\Lambda$ be the homoloidal net associated to $\varphi$ (which is the strict pull-back by $\varphi$ of the linear system of lines of $\Pn$) and let $\eta:X\rightarrow \Pn$ be a minimal birational morphism that solves the indeterminacies of $\varphi$ (or equivalently the base-points of $\Lambda$). 
Denote by $d$ the degree of $\varphi$ (which is the degree of the curves of $\Lambda$), by $q_1,...,q_k$ the base-points of $\varphi$ (or $\Lambda$), that may be proper or infinitely near points of $\Pn$, and by $a_i$ the multiplicity of $q_i$ as a base-point of $\Lambda$. We have $a_i\geq 1$ and $m_{q_i}\geq 0$.

Consider now  the strict transforms $\widetilde{\Lambda}$ of $\Lambda$  and $\widetilde{C}$ of $C$ on $X$. Then, $\widetilde{\Lambda}$ is base-point-free and $\widetilde{\Lambda}^2=1$. Using adjunction formula we find the classical equality $3(d-1)=\sum_{i=1}^k a_i$. Computing the free intersection of $\widetilde{\Lambda}$ and $\widetilde{C}$ (which is equal to the degree of the image $D$ of $C$, and is, by the hypothesis, at most equal to $n$), we find $dn-\sum_{i=1}^k a_i\cdot m_{q_i}\leq n$. This yields, with the equality above:
\begin{equation}\sum_{i=1}^k n\cdot a_i=3n(d-1)\leq\sum_{i=1}^k 3m_{q_i} \cdot a_i.\label{inequationT}\end{equation}
Since $3m_{q_i}\leq n$ and $a_i\geq 1$ for $i=1,...,k$, the inequality (\ref{inequationT}) is an equality. This implies that $\deg(D)=n$ and  $3m_{q_i}=n$ for $i=1,..,k$.
\end{proof}

\begin{cor}\label{Cor:SmoothCurves}
Let $C\subset \Pn$ be a smooth curve of degree $n$.
\begin{enumerate}
\item
If $n=3$, every base-point of every element of $\dec(C)$ belongs to $C$, as a proper or infinitely near point. 
\item
If $n>3$, then every element of $\dec(C)$ is an automorphism of the plane, i.e. $\dec(C)=\Aut(\Pn,C)$.\end{enumerate}
\end{cor}
\begin{proof} Apply Theorem~\ref{Thm:BasePoints}, with $m_q=1$ for any point $q$ that belongs to $C$ as a proper or infinitely near point.\end{proof}
The first part of Corollary~\ref{Cor:SmoothCurves} may be found in \cite[Theorem 1.3]{Pan} and the second in \cite[Cor. 3.6]{Pan} and \cite{Kup} (see also \cite[p. 181]{Rep} and \cite[Book IV, Chap. VII, \S 3, Thm. 11]{Coo}).

Another important Corollary is the following one, that describes the inertia group of a family of classical curves (Halphen curves, Coble curves, ...) as a subgroup of automorphisms of a rational surface. We will use this to provide examples of plane curves whose canonical complex is not exact.

\begin{cor}\label{Cor:Aut}
Let  $p_1,...,p_k\in\Pn$ be $k$ distinct proper points of $\Pn$ and let $C\subset \Pn$ be an irreducible curve of degree $3n$, with $n>1$, and which has multiplicity $n$ at each $p_i$. Denote by $\pi:X\rightarrow \Pn$ the blow-up of the $k$ points and assume that the strict pull-back $\widetilde{C}$ of $C$ by $\pi$ is a smooth curve. 
Then, 

\begin{enumerate}
\item
$\pi^{-1} \dec(\Pn,C)\pi=\dec(X,\widetilde{C})=\Aut(X,\tilde{C})$.

\item Let $D\subset \Pn$ be an irreducible curve of degree $\leq 3n$ and $\varphi:(\Pn,C)\tor(\Pn,D)$ a birational map. 
Denote by $\eta:Y\rightarrow \Pn$ an embedded minimal resolution of the singularities of $D$ and by $\tilde{D}\subset Y$ the strict transform of $D$. 

Then, $\varphi$ lifts to an isomorphism $\varphi':(X,\widetilde{C})\to (Y,\widetilde{D})$ such that $\eta\varphi'=\varphi\pi$. Furthermore, the degree of $D$ is $3n$.
\end{enumerate}
\end{cor}
\begin{proof}
Let us prove assertion $(2)$. Theorem~\ref{Thm:BasePoints} implies that the base-point locus of $\varphi$ is contained in $\{p_1,...,p_k\}$. Then, $\varphi\pi$ is a birational morphism $X\rightarrow \Pn$, that we denote by $\nu$. Since the curve $\tilde{C}$ is equivalent to $-nK_X$, the degree of the curve $\nu(\tilde{C})=D$ is $3n$ and every $(-1)$-curve of $X$ intersects $\tilde{C}$ into $n$ points. This implies that $\nu$ is an embedded minimal resolution of the singularities of $D$. The two birational morphisms $\nu$ and $\eta$ differ only by an isomorphism $\varphi':X\rightarrow Y$, which sends $\tilde{C}$ on $\tilde{D}$.

The assertion $(1)$ is a particular case of (2). Indeed, for $\phi\in\dec(\Pn,C)$, the element  $\phi'=\pi^{-1}\phi\pi$ belongs to $\Aut(X,\tilde{C})$ and consequently the group $\pi^{-1}\dec(\Pn,C)\pi$ is contained in $\Aut(X,\tilde{C})$; the other inclusion is obvious.
\end{proof}

\subsection{Examples of different birational embeddings of curves of genus $2$ in $\Pn$}
Let $C$ be any abstract smooth curve of genus $2$. It is isomorphic to the curve $z^2=F_6(x,y)$ in the weighted projective plane $\mathbb{P}(1,1,3)$, for some form $F_6$ of degree $6$, having $6$ simple roots. There exists a birational morphism $C\stackrel{\varphi_1}{\rightarrow} C_0$ where $C_0$ is a quartic of $\Pn$ with one ordinary double point, and furthermore there is only one choice of $C_0$, up to birational equivalence of the pair $(\Pn,C_0)$ (see \cite{BayBea}). The group $\ine(C_0)$ is infinite and described in Example~\ref{Exa:deJAbGr} (Theorem~\ref{THM:BPVClass}); moreover the canonical complex of $(\Pn,C_0)$ is an exact sequence (Theorem~\ref{THM:SplitExactSequenceWhenIne1}).

Let $F_4$ be any form of degree $4$ in two variables (that may also be equal to zero), and define $S$ to be the surface with equation $w^2=z^3+zF_4(x,y)+F_6(x,y)$ in the weighted projective space $\mathbb{P}(3,1,1,2)$. Since $F_6$ does not have any multiple root, $S$ is smooth and then is a del Pezzo surface of degree~$1$ \cite[Theorem 3.36]{KoSmCo}. There exists thus a birational morphism $\pi:S\rightarrow \Pn$ that consists of blowing-up $8$ points in general position. Sending the curve $C$ into $S$ via the morphism $(x:y:z)\mapsto (z:x:y:0)$ gives a curve $\tilde{C}\subset S$, equivalent to $-2K_S$, whose image by $\pi$ is a sextic with eight ordinary double points. 

Assume that $F_4$ is the zero form and denote by $C_1\subset \Pn$ the image of $\tilde{C}$ by $\pi$. Then $\ine(C_1)$ is isomorphic to $\Z/3\Z$, and the canonical complex of $(\Pn,C_1)$ is an exact sequence (Theorems~\ref{THM:BPVClass} and~\ref{THM:SplitExactSequenceWhenIne1}).

Assume that $F_4$ is not the zero form, then no non-trivial automorphism of $S$ fixes the curve $\tilde{C}$, which means that $\ine(S,\widetilde{C})$ is trivial (Corollary~\ref{Cor:Aut}). The Bertini involution on $S$ (that sends $w$ on $-w$) leaves $\tilde{C}$ invariant, acts on it as the involution associated to the $g^1_2$, and generates $\dec(S,\tilde{C})$ if $F_4$ is general enough. Moreover $\Aut(\tilde{C})$ is reduced to this involution if and only if no non-trivial automorphism of $\mathbb{P}^1$ leaves $F_6$ invariant. It follows then from Corollary~\ref{Cor:Aut} that the canonical complex is an exact sequence under these circumstances, here $C_2$ denotes the image of $\tilde{C}$ by $\pi$.

 Theses examples provide three birational embeddings $C\rightarrow C_i\subset  \Pn$ that leads to three different canonical complexes for the same abstract curve and also to  three birationally different pairs $(\Pn,C_i)$. Theorem~\ref{Thm:BasePoints} allows us to improve this result, giving infinitely many such pairs of the last kind. Indeed, let $F_4$ and $F_4'$ be two different forms of degree~$4$, let $\tilde{C}\subset S$ and $\tilde{C}'\subset S'$ be the two embeddings of $C$ into two corresponding del Pezzo surfaces of degree~$1$, and let $C_2\subset \Pn$ and $C_2'\subset \Pn$ be the corresponding sextic curves. If there exists a birational transformation $\varphi$ that sends $(\Pn,C_2)$ on $(\Pn,C_2')$, then Corollary~\ref{Cor:Aut} implies that $\varphi$ lifts to an isomorphism $S\rightarrow S'$. Changing the values of $F_4$, we reach infinitely many isomorphism classes of del Pezzo surfaces of degree $1$, that leads to infinitely many birationally different pairs $(\Pn,C_2)$ such that $C_2$ is birational to $C$.

\subsection{Examples of different birational embeddings of curves of genus $3$ in $\Pn$}\label{SubSecGenus3}
We give another example. Let $C_1\subset \Pn$ be any smooth quartic curve. The double covering of $\Pn$ ramified over $C_1$ is a del Pezzo surface $S$ of degree $2$ (see \cite{bib:BeauLivre}), which is the blow-up $\pi:S\rightarrow \Pn$ of $7$ points of $\Pn$ in general position. Denote by $\tilde{C}$ the image of $C_1$ on $S$ and by $C_2$ the curve $\pi(\tilde{C})$. Then, $C_2$ is a sextic with $7$ ordinary double points and $\ine(C_2)\cong \Z/2\Z$ is generated by the Geiser involution that corresponds to the involution of $S$ associated to the double covering (Theorem~\ref{THM:BPVClass}). On the other hand, Corollary~\ref{Cor:SmoothCurves} implies that $\dec(C_1)=\Aut(\Pn,C_1)$ and consequently that $\ine(C_1)$ is trivial. The two curves $C_1$ and $C_2$ are birational curves of the plane, but the pairs $(\Pn,C_1)$ and $(\Pn,C_2)$ have different canonical complexes, and in particular are not birationally equivalent.

\subsection{Examples of different birational embeddings of curves of genus $4$ in $\Pn$}

Let $p_1,...,p_8$ be eight points of the plane and let $S\rightarrow \Pn$ be the blow-up of these points. Assume that $S$ is a del Pezzo surface. Corollary~\ref{Cor:Aut} implies the following observations. Among the linear system $\Lambda$ of nonics passing through $p_1,...,p_8$ with multiplicity 3, exactly one has a non-trivial inertia group, generated by the Bertini involution of $S$. The other curves of $\Lambda$ have a decomposition group that contains the Bertini involution, and for a general curve of $\Lambda$, this involution generates the decomposition group. Furthermore, the elements of $\Lambda$ yields infinitely many pairs which are birationally different.

\section{Curves of genus one}\label{Sec:Genus1}
In section~\ref{Sec:GenusAtLeast2}, we gave a precise description of all elements of finite order of $\cre(\Pn)$ that fix a curve of genus $\geq 2$. A so precise description exists for curves of genus $1$:
\begin{thm}\label{theo:FiniteOrderGenus1}\PrpTitle{\cite[Theorem 2]{BlaMich}}

Let $C\subset \Pn$ be an irreducible curve with $g(C)=1$. Let $h\in \ine(C)$ be an element of finite order $n>1$. Then, there exists a birational map $\varphi:\Pn\dasharrow S$ that conjugates $h$ to an automorphism $\alpha$ of a del Pezzo surface $S$, such that $(\alpha,S)$ are given in the following table:

\begin{tabular}{|p{0.1 mm}p{3 mm}p{2 cm}p{6.8 cm}p{2 cm}|}
\hline
&n\parbox{3 mm}{{\color{white}n \hts}\\ {\it  \hts}}& \parbox{25 mm}{{\it description}\\ {\it of $\mathit{\alpha}$}} & \parbox{62 mm}{{\it equation of}\\ {\it the surface $S$}}& \parbox{62 mm}{{\it in the}\\ {\it variety}} \\
\hline
&2& $x_0\mapsto -x_0$& $\sum_{i=0}^4 x_i^2=\sum_{i=0}^4 \lambda_i x_i^2=0$& $\mathbb{P}^4$  \\
\hline 
&3& $x_0\mapsto \zeta_3x_0$& ${x_0}^3+L_3(x_1,x_2,x_3)$& $\mathbb{P}^3$  \\
\hline 
&4& $x_0\mapsto \zeta_4 x_0$ & ${x_3}^2={x_0}^4+L_4(x_1,x_2)$ & $\mathbb{P}(1,1,1,2)$  \\
\hline 
&5& $x_0\mapsto \zeta_5 x_0$ & ${x_3}^2={x_2}^3+\lambda_1 {x_1}^4x_2+{x_1}(\lambda_2 {x_1}^5+{x_0}^5)$ & $\mathbb{P}(1,1,2,3)$\\
\hline 
&6& $x_0\mapsto \zeta_6 x_0$ & ${x_3}^2={x_2}^3+\lambda_1 {x_1}^4x_2+\lambda_2 {x_1}^6+{x_0}^6$ & $\mathbb{P}(1,1,2,3)$,  \\
 \hline
\end{tabular}
where $\zeta_n\in \C$ is a primitive $n$-th root of the unity, $L_i$ is a form of degree $i$ and $\lambda_i$ are parameters such that $S$ is smooth. 

Furthermore, any birational morphism $S\rightarrow \Pn$ sends the fixed curve on a smooth plane cubic curve.
\end{thm}
Theorem~\ref{theo:FiniteOrderGenus1} implies in particular the following result:
\begin{cor}\label{Cor:IneCubic23456}
Let $C\subset \plan$ be an irreducible curve with $g(C)=1$. The following conditions are equivalent:
\begin{enumerate}
\item
the pair $(\Pn,C)$ is birationally equivalent to a pair $(\Pn,D)$, where $D$ is a smooth cubic curve.
\item
the group $\ine(C)$ contains elements of finite order;
\item
the group $\ine(C)$ contains elements of order $2$, $3$, $4$, $5$ and $6$.
\end{enumerate}
\end{cor}

\begin{proof}
In order to prove $(1)\Rightarrow(3)$ we observe that in each of the five types of pairs $(\alpha, S)$ of Theorem~\ref{theo:FiniteOrderGenus1}, one gets an arbitrary elliptic curve as fixed curve. The implication $(3)\Rightarrow(2)$ is obvious and $(2)\Rightarrow(1)$ follows from Theorem~\ref{theo:FiniteOrderGenus1}.
\end{proof}

The curves of genus $1$ having the biggest canonical complex seems in fact to be the cubic curves. We will precise it in (\ref{SubSec:Genus1NotCubic}). We examine in (\ref{SubSec:SmoothCubIne}) and (\ref{SubSec:SmoothCubDec}) the case of smooth cubic curves and then in (\ref{SubSec:Genus1NotCubic}) the other irreducible curves of genus $1$. 

\subsection{The inertia group of smooth cubic curves}\label{SubSec:SmoothCubIne}
Let $C\subset \Pn$ be a smooth cubic curve. Taking any point $p\in C$, there exist infinitely many elements of $\ine(C)$ that leave invariant any general line passing through $p$; such elements form a group which is described in Example~\ref{Exa:deJAbGr}. There are furthermore many elements of degree $3$ in this group (\cite[Lemma 4.1]{BlaMich}), one of these, that we call $\sigma_p$ is the classical de Jonqui\`eres involution of Example~\ref{Exa:GeiBerJon}c) (generalised in \cite{Giz} over the name of $R_p$ to any dimension).  The element $\sigma_p$ is the unique involution that leaves invariant any general line passing through $p$ and fixes the curve $C$. 

Changing the choice of $p$, all these involutions generate a very large group:
\begin{thm}\PrpTitle{\cite[Theorem 1.6]{BlaMich}}\label{ThmCubicFree}
Let $C\subset \Pn$ be a smooth cubic curve. The subgroup of $\ine(C)$ generated by all the cubic involutions centred at the points of $C$ is the free product $${{\star}}_{p\in C} <\sigma_p>.$$
\end{thm}

Furthermore, since the inertia group of a smooth cubic curve contains elements of order $3$, $4$, $5$ and $6$ (Corollary~\ref{Cor:IneCubic23456}), the free product described in Theorem~\ref{ThmCubicFree} is not the whole inertia group.
However, there exists an analogous of Noether-Castelnuovo's theorem for this group:
\begin{thm}\PrpTitle{\cite[Theorem 1.1]{BlaMich}}\label{ThmCubicCastelnuovo}
The inertia group of a smooth plane cubic curve is generated by its elements of degree $3$, which are -- except the identity -- its elements of lower degree.
\end{thm}

%


\subsection{The decomposition group of smooth cubic curves}\label{SubSec:SmoothCubDec}
Let $C$ be a smooth plane cubic curve. Take three distincts points $p,q,r$  that belong to $C$ as proper or infinitely near points. The linear system of conics passing through these points defines a birational transformation $\varphi$ of $\plan$ which transforms $C$ onto another smooth cubic curve $C'$. Composing $\varphi$ with a linear automorphism mapping $C'$ onto $C$ we obtain a degree 2 element in $\dec(C)$. Clearly these transformations are the only degree 2 elements in $\dec(C)$. Moreover, all such transformations may be expressed as composition of those whose base-point set consists of three proper points of the plane. Like for the inertia group, there exists an analogous of Noether-Castelnuovo's theorem for the decomposition group:
\begin{thm}\PrpTitle{\cite[Theorem 1.4]{Pan}}\label{ThmCubicDecom}
The decomposition  group of a smooth plane cubic curve is generated by its elements of degree $2$.
\end{thm}
Concerning the action of the decomposition group on the elliptic curve, the following shows that this one is like the whole automorphism group of the curve:
\begin{thm}\PrpTitle{\cite[Theorem 6]{Giz}}\label{ThmCubicComplex}
Let $C\subset \Pn$ be a smooth cubic curve. The canonical complex of $(\Pn,C)$ is an exact sequence.
\end{thm}
\begin{rem}
It seems that the sequence is not split.
\end{rem}

\subsection{Curves of genus $1$ that are not equivalent to smooth cubic curves}\label{SubSec:Genus1NotCubic}
Recall some classical notions on Halphen curves and surfaces (see \cite{Hal}, \cite{Cob}, \cite{GizH}, \cite{DolOrt}).
\begin{defi}
A \emph{Halphen curve of index $n$} is an irreducible plane curve of degree $3n$, with $9$ points of multiplicity $n$ and of genus $1$. 

A projective rational smooth surface $S$ is a \emph{Halphen surface  of degree $n$} if the linear system $|-nK_S|$ is a pencil whose general fibre is an irreducible curve of genus $1$.
\end{defi}

The following classical relation may be verified by hand:
\begin{lem}
If $S$ is a Halphen surface of degree $n$, any birational morphism $S\rightarrow \Pn$ sends the general fiber of $|-nK_S|$ on Halphen curves of index $n$.

For $n\geq 2$, the blow-up of the $9$ singular points of a Halphen curve of index $n$ is a Halphen surface.
\end{lem}
The blow-up of $9$ general points is not a Halphen surface. However, for any general set of $8$ points of the plane, and for any integer $n\geq 2$, there exists a curve of the plane such that the blow-up of the $8$ points and a ninth point on the curve gives a Halphen surface of index $n$ \cite{Hal}.

We give now a simple proof to the following result, that is probably classical.
\begin{pro}
Let $C_1,C_2\subset \Pn$ be two Halphen curves of index respectively $n_1$ and $n_2$. For $i=1,2$, let $\eta_1:X_1\rightarrow \Pn$ be the minimal embedded resolution of $C_i$ (which is the identity if $n_i=1$). The following assertions are equivalent:
\begin{enumerate}
\item
the pairs $(\Pn,C_1)$ and $(\Pn,C_2)$ are birationally equivalent;
\item
there exists an isomorphism $\varphi:X_1\rightarrow X_2$ that sends the strict transform of $C_1$ on the strict transform of $C_2$.
\end{enumerate}
Furthermore, both assertions imply that $n_1=n_2$.
\end{pro}
\begin{proof}
The second assertion implies directly the first one and the equality $n_1=n_2$. Corollary~\ref{Cor:Aut} shows that the first assertion implies the second one.
\end{proof}
This proposition shows in particular the existence of infinitely many distinct types of pairs $(\Pn,C)$ where $C$  has genus $1$; it raises also the following question, which seems to us open until now.
\begin{ques}
Let $C_1\subset \Pn$ be an irreducible curve of genus $1$. Does there exists a Halphen curve $C\subset \Pn$ such that the pair $(\Pn,C_1)$ is birationally equivalent to $(\Pn,C)$?
\end{ques}

We describe now the decomposition and inertia groups of Halphen curves of index $\geq 2$ (those of index $1$ are the smooth cubic curves, described previously), and show in particular the important difference between index $1$ and index $\geq 2$.
\begin{thm}\label{THM:HalphZ8}
Let $C\subset \Pn$ be a Halphen curve of index $n\geq 2$. Then, $\dec(C)$ contains a normal subgroup of finite index, isomorphic to $\mathbb{Z}^8$.
In particular, the canonical complex of $(\Pn,C)$ is not exact.

Assume that $C$ is a general Halphen curve; then either $\dec(C)$ is isomorphic to $\mathbb{Z}^8\rtimes \mathbb{Z}/2\mathbb{Z}$ or to $\mathbb{Z}^8$. The first case occurs for  $n=2$ and the second one  if $n=3$ and $n\geq 5$.
\end{thm}
\begin{proof}
Let $\pi:S\rightarrow \Pn$ be the blow-up of the nine singular points of $C$ and let $\tilde{C}\subset S$ be the strict transform of $C$. Corollary~\ref{Cor:Aut} implies that $\dec(C)=\dec(\Pn,C)$ is conjugate to $\dec(S,\tilde{C})=\Aut(S,\tilde{C})$. 

Denote by $D\subset \Pn$ a cubic passing through the singular points of $C$ and by $\tilde{D}\subset S$ the strict pull-back of $D$. Then $\tilde{C}$ and $n\tilde{D}$ take part of the pencil $|-nK_S|$ and this shows that $D$ is unique.

Note that $\Aut(S)$ acts on the elliptic fibration $\eta:S\rightarrow \mathbb{P}^1$ induced by $|-nK_S|$. Let $G\subset \Aut(S)$ be the subgroup of automorphisms that act trivially on the basis and let $G'$ be the image of $\Aut(S)$ in $\Aut(\mathbb{P}^1)$, such that the following is an exact sequence:
\[1\rightarrow G\rightarrow \Aut(S)\rightarrow G'\rightarrow 1.\]

We show that $G'$ is finite. Indeed, $\eta(\tilde{D})$ is a fixed point, so we can consider $G'$ as a subgroup of $\Aut(\C)$; then $G'$ has at least one finite orbit in $\C$ because there are singular fibres in $|-nK_S|$ (the Euler characteristic of $S$, which is equal to $12$, is the sum of the Euler characteristics of the (singular) fibres of $\eta$), which is not possible if $G'$ is infinite.

Now, let $H\subset G$ be the subgroup of elements that act as translations on the general fibre; according to the structure of the automorphisms group of an elliptic curve and since $G'$ is finite, $H$ is normal in $\Aut(S)$, of finite index.

A translation on an elliptic curve corresponds to a linear equivalence of a divisor of degree $0$. There exists thus an exact sequence (see \cite{GizH})
\[0\rightarrow \mathbb{Z}K_S\rightarrow K_S^{\perp}\rightarrow H\rightarrow 0,\]
where $K_S^{\perp}$ is the subgroup of $\pic(S)$ of elements whose intersection with $K_S$ is equal to $0$. Since $\pic(S)\cong \mathbb{Z}^{10}$ and $K_S$ is indivisible, ${K_S}^{\perp}\cong \mathbb{Z}^9$ and $H\cong \mathbb{Z}^8$.
As $H\subset G\subset \dec(S,\tilde{C})\subset \Aut(S)$, the first assertion is proved.

Assume now that $C$ is a general Halphen curve, which implies that $\dec(S,\tilde{C})$ is equal to $G$ and that the automorphism group of $C$ is equal to $C\rtimes \mathbb{Z}/2\mathbb{Z}$. In particular, the index of $H$ in $G$ is either $1$ or $2$, depending on whether there exists an element of $G$ that acts as an involution with four fixed points on the general fibre. Assume that such an element $\sigma\in G$ exists. Then, it fixes a curve in $S$ which intersects the general fibre into four points. Since the fibre is equal to $-nK_S$, then $n$ must divide $4$, which implies that $n=2$ or $n=4$. 

It remains to show that for $n=2$, such an involution exists. 
 Consider the elliptic fibration $\epsilon:S\dasharrow \mathbb{P}^1$ defined by the pencil of plane cubics passing through eight of the nine base-points of the Halphen pencil: The intersection of a general fiber $S_{\eta}$ of  $\eta$ with a general fiber of $\epsilon$ is equal to two, which means that the degree of $\eta \times \epsilon:S\dasharrow \mathbb{P}^1\times\mathbb{P}^1$ is two; the corresponding  
involution of $S$ leaves each Halphen curve $S_{\eta}$ invariant and has thereon four fixed points since the restriction of  $S_{\eta}\dasharrow \mathbb{P}^1$ of $\epsilon$ is nothing but the canonical $g_1^2$.
\end{proof}
\begin{rem}
In the case where $n=2$ and the points are in general position, it may also be observed that the Bertini involution associated to the blow-up of $8$ of the $9$ points lifts to an automorphism of the surface that acts on each member of the elliptic fibration as an automorphism with four fixed points. Furthermore, the $9$ involutions obtained via this map generate the automorphism group of the Halphen surface \cite{Cob}.

\end{rem}
\section{Rational curves} \label{Sec:Genus0}
The case of rational curves is the less described. We can however cite some simple results.
\subsection{The line}
There exists plenty of elements in the inertia group of the line: for example, any birational map of the form $(x,y)\dasharrow \left(\frac{x}{\alpha(y)x+\beta(y)},y\right)$, where $\alpha,\beta \in \C(x), \beta\not=0$, fixes the line $x=0$. It seems that the inertia group of a line is a big and complicated group. Let us give some simple observations:
\begin{pro}
Let $L\subset \Pn$ be a line, then the canonical complex of $(\Pn,L)$ is a split exact sequence.

Furthermore, the group $\ine(\Pn,L)$ is neither finite, nor abelian  and does not leave invariant any pencil of rational curves.
\end{pro}\begin{proof}
The exactness and splitness are obvious: the group of automorphisms of $L$ extends to a subgroup of $\Aut(\Pn)$, and this yields a section $\Aut(L)\rightarrow \Aut(\Pn,L)$.
 The other assertions follow from \cite[Proposition 4.1]{BPV}.
\end{proof}
Does there exists an analogous of Noether-Castelnuovo's theorem, as for the case of smooth cubics?
 \begin{ques}
Let $L\subset \Pn$ be a line. Is the group $\dec(\Pn,L)$ (respectively $\ine(\Pn,L)$) generated by its elements of degree $1$ and $2$?
\end{ques}

\subsection{Coble curves}\label{SubSecCoble}
A Coble curve is an irreducible sextic with $10$ double points. There does not exist a cubic singular at ten general points; however Coble curves exist, and are singular members of a Halphen pencil of index $2$; furthermore in each such pencil there are in general $12$ Coble curves \cite{Hal}. Corollary~\ref{Cor:Aut} implies that the pair $(\Pn,C)$ where $C$ is a Coble curve is not equivalent to the one of a line. Furthermore, we have the following:
\begin{pro}
Let $C\subset \Pn$ be a Coble curve, let $\pi:S\rightarrow \Pn$ be the blow-up of its $10$ singular points and let $\tilde{C}\subset S$ be the strict pull-back of $C$ by $\pi$.

Then, $\Aut(S)=\dec(S,\tilde{C})=\pi^{-1}\dec(\Pn,C)\pi$.
\end{pro}
\begin{proof}The curve $\tilde{C}$ is equivalent to $-2K_S$ and since it has negative self-intersection it is the only such curve, consequently, 
$\Aut(S)=\Aut(S,\tilde{C})$. The result follows then directly from Corollary~\ref{Cor:Aut}.
\end{proof}

The description of the automorphisms of a so-called Coble surface obtained by blowing-up the ten singular points of a Coble curve is a classical result of Coble \cite{Cob}, see also \cite[Theorem 8, page 107]{DolOrt}. It implies in particular the following result:
\begin{pro}
For any Coble curve $C$, the group $\dec(\Pn,C)$ is an infinite countable group. The canonical complex associated to $(\Pn,C)$ is not exact.
\end{pro}

\subsection{Other curves of Halphen type}\label{SubSecRationalofHalphenType}
Let $S$ be a Halphen surface of index $n$ obtained by the blow-up $\pi:S\rightarrow \Pn$ of the points
$p_1,...,p_{9}$. There exist singular fibres of the elliptic fibration $|-nK_S|$, which are thus rational curves
with a double point, whose image on $\Pn$ are curves of degree $3n$ with multiplicity $n$ at the points
$p_1,...,p_9$ and multiplicity $2$ at some other point $p_{10}$. The case $n=1$ gives nodal cubics, which are
equivalent to lines; the case $n=2$ gives Coble curves, and the case $n=3$ gives other curves. Once again, it seems
that in general $12$ such curves exist in a general Halphen pencil \cite{Hal}.

\begin{pro}
Let $C\subset \Pn$ be an irreducible curve of degree $3n$ with multiplicity $n$ at $p_1,...,p_9$ and multiplicity $2$ at $p_{10}$, and assume that $n\geq 3$.
Let $\pi:S\rightarrow \Pn$ (respectively $\pi':S'\rightarrow \Pn$) be the blow-up of $p_1,...,p_{10}$ (respectively of $p_1,...,p_{9}$), and let $\tilde{C}\subset S$ and $\tilde{C}'\subset S'$ be the strict pull-backs of $C$ by $\pi$ and $\pi'$.

Then, $\Aut(S)=\dec(S,\tilde{C})=\pi^{-1}\dec(\Pn,C)\pi$ and $\Aut(S',\tilde{C}')=\dec(S,\tilde{C'})=\pi'^{-1}\dec(\Pn,C)\pi'$.

Furthermore, $\dec(\Pn,C)$ contains a subgroup of finite index isomorphic to $\mathbb{Z}^8$. In particular, the canonical complex associated to $(\Pn,C)$ is not an exact sequence.
\end{pro}
\begin{proof}
As for Corollary~\ref{Cor:Aut}, Theorem~\ref{Thm:BasePoints}  implies the equalities  $\Aut(S,\tilde{C})=\dec(S,\tilde{C})=\pi^{-1}\dec(\Pn,C)\pi$. The curve $E_{10}=\pi^{-1}(p_{10})$ is a smooth irreducible rational curve of self-intersection $-1$ (a $(-1)$-curve) and its intersection with $\tilde{C}$ is $2$; furthermore, it is the unique such curve \cite[Theorem~3.3]{KuM}. Consequently, the whole group $\Aut(S)$ leaves $E_{10}$ invariant; denoting by $\mu:S\rightarrow S'$ the blow-down of this curve (such that $\pi=\pi'\mu$), the group $G=\mu \Aut(S)\mu^{-1}$ is the subgroup of $\Aut(S')$ that fixes the point $(\pi')^{-1}(p_{10})=\mu(E_{10})$, which is the unique singular point of $\tilde{C}'$. Since $S'$ is a Halphen surface of index $n$ and 
$\tilde{C}'$ is a singular member of the fibration, $G=\Aut(S',\tilde{C}')$. This implies the remaining equalities.

The last part follows from Theorem~\ref{THM:HalphZ8}.
\end{proof}
\subsection{Other rational curves}
Do there exist other examples of pairs $(\Pn,C)$ where $C$ is rational? A famous Coolidge-Nagata problem asks whether the pair of a rational cuspidal curve is birationally equivalent to the pair of a line (see \cite{Coo} and \cite{Nag}). 

\begin{defi}
Let $C\subset S$ be an irreducible smooth curve in a surface. We denote by $\kappa(S,C)$ the Kodaira dimension of the pair $(S,C)$: this is the dimension of the image of $X\dasharrow \mathbb{P}(H^0(m(D+K))^{\vee})$ for $m$ large enough. If $|m(C+K_S)|=\emptyset$ for any $m>0$, the Kodaira dimension is by convention equal to $-\infty$.

For a singular curve $C\subset S$, we write $\kappa(S,C)=\kappa(X,\tilde{C})$, where $X\rightarrow S$ is an embedded resolution of the singularities of $C$ and $\tilde{C}\subset X$ is the strict transform.
\end{defi}

\begin{lem}\PrpTitle{\cite{KuM}}
If $(S,C)$ is birationally equivalent to $(S',C')$ then $\kappa(S,C)=\kappa(S',C')$.
\end{lem}

Let us cite the following fondamental result, due to Coolidge. 
\begin{thm}\PrpTitle{\cite{Coo}, \cite{KuM}}
Let $C\subset \Pn$ be a rational irreducible curve and $L\subset \Pn$ be a line. Then $(\Pn,C)$ is birationally equivalent to $(\Pn,L)$ if and only if $\kappa(\Pn,C)=-\infty$.
\end{thm}

We have also a description for Kodaira dimension $0$ and $1$:
\begin{thm}\PrpTitle{\cite{KuM}}
Let $C\subset \Pn$ be a rational irreducible curve.

\begin{enumerate}
\item
$\kappa(\Pn,C)=0$ if and only if $(\Pn,C)$ is birationally equivalent to $(\Pn,D)$ where $D$ is a Coble curve.
\item
$\kappa(\Pn,C)=1$ if and only if $(\Pn,C)$ is birationally equivalent to $(\Pn,D)$, where $D$ is a curve of degree $3n$, with $9$ points of multiplicity $n>2$ and a tenth point of multiplicity $2$.
\end{enumerate}
\end{thm}

Consequently, finding other rational curves not equivalent to our examples is equivalent to find rational curves $C\subset \Pn$ with $\kappa(\Pn,C)=2$.
\section{Link between the inertia and decomposition groups and the classification of finite subgroups of the Cremona group}\label{Sec:LinkFinite}
In our description of the decomposition group, and more precisely the inertia group of a curve of genus $\geq 1$, we provide some subgroups of finite order. Conversely, in the study of the finite subgroups of $\Bir(\Pn)$, the set of birational classes of curves of positive genus fixed (pointwise) is an important conjugacy invariant.

\subsection{Cyclic groups of prime order}The conjugacy class of a finite cyclic subgroups of prime order of $\Bir(\Pn)$ is uniquely determined by the birational equivalence of the curve of positive genus that it fixes (it may fix at most one such a curve):

\begin{thm}\PrpTitle{\cite{BayBea}, \cite{Fer}}
Let $g,g'$ be two elements of $\Bir(\Pn)$ of the same prime order, that fix respectively $\Gamma, \Gamma'$, two irreducibles curves of positive genus.
Then, $g$ and $g'$ are conjugate in $\Bir(\Pn)$ if and only if the curves $\Gamma$ and $\Gamma'$ are birational. 
\end{thm}
\begin{thm}\PrpTitle{\cite{BeaBla}}\label{Thm:CurveFixed}
An element of $\Bir(\Pn)$ of prime order is not conjugate to a linear automorphism if and only if it belongs to the inertia group of some curve of positive genus. 
\end{thm}
\subsection{Other groups}
Theorem~\ref{Thm:CurveFixed} extends to finite cyclic group of any order, and almost to finite abelian groups:

\begin{thm}\PrpTitle{\cite{BlaLin}, announced in \cite{BlaCR}}
\label{Thm:Cyclic}
Let $G$ be a finite cyclic subgroup of order $n$ of $\Bir(\Pn)$. The following conditions are equivalent:
\begin{itemize}
\item
If $g \in G$, $g\not=1$, then $g$ does not fix a curve of positive genus.
\item
$G$ is birationally conjugate to a subgroup of $\Aut(\Pn)$.
\item
$G$ is birationally conjugate to a subgroup of $\Aut(\mathbb{P}^1\times\mathbb{P}^1)$.
\item
$G$ is birationally conjugate to the group of automorphisms of $\Pn$ generated by $(x:y:z) \mapsto (x:y:e^{2\im\pi/n} z)$.
\end{itemize}
\end{thm}
\begin{thm}\PrpTitle{\cite{BlaLin}, announced in \cite{BlaCR}}
\label{Thm:NonCyclic}
Let $G$ be a finite abelian subgroup of $\Bir(\Pn)$. The following conditions are equivalent:
\begin{itemize}
\item
If $g \in G$, $g\not=1$, then $g$ does not fix a curve of positive genus.
\item
$G$ is birationally conjugate to a subgroup of $\Aut(\Pn)$, or to a subgroup of $\Aut(\mathbb{P}^1\times\mathbb{P}^1)$ or to the group $\mathit{Cs}_{24}$ isomorphic to $\mathbb{Z}/2\mathbb{Z}\times\mathbb{Z}/4\mathbb{Z}$, generated by the two elements \[\begin{array}{lll}(x:y:z)&\dasharrow&(yz:xy:-xz),\\
(x:y:z)&\dasharrow &( yz(y-z):xz(y+z):xy(y+z)).\end{array}\]
\end{itemize}
Moreover, this last group is conjugate neither to a subgroup of $\Aut(\Pn)$, nor to a subgroup of $\Aut(\mathbb{P}^1\times\mathbb{P}^1)$.
\end{thm}

However, there are plenty of examples of finite non-abelian subgroups of $\Bir(\Pn)$ which are not birationally conjugate to a subgroup of $\Aut(\Pn)$ or $\Aut(\mathbb{P}^1\times\mathbb{P}^1)$ but such that no curve of positive genus is fixed by any non-trivial element of the group \cite[Section 12]{BlaLin}.

\section{The links with the dynamical properties of a Cremona transformation}\label{Sec:Dynamic}
We can also consider a Cremona transformation as defining a dynamical system. In comparison with the usual case of dynamics defined by automorphisms, the situation here is more complicated due to the presence of indeterminacies and critical points: in the neighbourhood of such points the map does not act in a "natural way".

In \cite{Fri} and \cite{RuSh} the authors introduce the so-called first dynamical degree of a birational map; this number is invariant by birational conjugation. Let us explain what that degree means in the case of a Cremona transformation $\varphi$. Consider the sequence $\sqrt[n]{\deg(\varphi^n)}$ for $n\geq 0$. Since $\deg(\varphi^{n+m})\leq \deg(\varphi^n) \cdot\deg(\varphi^{m})$ it has a limit. The first dynamical degree is then
\[\lambda(\varphi):=\lim\,\sqrt[n]{\deg(\varphi^n)}.\]
As shown in \cite{Fri} the \emph{Topological Entropy} $h_{top}(\varphi)$ of $\varphi$ is at most $\log\,\lambda(\varphi)$. The equality is conjectured, and proved for a general $\varphi$ (\cite{BeDi} and \cite{Duj}).

On the other hand, Diller and Favre propose a more refined approach and consider the sequence of successive degrees $\deg(\varphi^n)$ itself. They classify the plane Cremona transformations (in fact they consider a more general setup) with $\lambda=1$ in terms of the growth rate of that sequence (see \cite[Thm. 0.2]{DiFa}); they show that this growth is at most quadratic in $n$.

It is natural to ask ourself in what extent the dynamic of a Cremona transformation may be affected by the existence of a genus $g$ stable curve, $g=0,1,2..$. One may answer the following:

\begin{thm}\PrpTitle{\cite[Theorem 1.1]{PanCR}}
Let $C\subset\plan$ be an irreducible curve of genus $g(C)\geq 2$ and let $\varphi\in\dec(C)$, then  $\lambda(\varphi)=1$ and the sequence $\{\deg(\varphi^n)\}_{n=1}^{\infty}$ grows at most linearly.
\end{thm}
\begin{proof}
Since $g(C)>1$, the subgroup $\ine(C)$ is of finite index in $\dec(C)$ so we may assume that $\varphi\in \ine(C)$. If $\varphi$ is of finite order, we are done. If not, Theorem~\ref{THM:BPVClass} yields an explicit description of $\varphi$ as in Example~\ref{Exa:deJAbGr}. Computing the degrees of the powers of $\varphi$, we have $\deg(\varphi^n)\leq n(\deg(\varphi)+c)$ for some constant $c$. This completes the proof.
\end{proof}
On the other hand, when $g(C)\leq 1$ the number $\lambda(\varphi)$ may be strictly larger than $1$.  For an example in the rational case we refer the reader to  \cite[Example 2]{PanCR}. In the case of a smooth cubic curve, the composition of two generic quadratic elements of the decomposition group seems to suit. 

Finally, until recently, all known examples of automorphisms of rational surfaces with a first dynamical degree strictly larger than $1$ (or equivalently with an action of infinite order on the Picard group) are those which leave invariant a rational or an elliptic curve. A question/ conjecture of Gizatullin/ Harbourne/ McMullen (see \cite{Har} and \cite{McM}) asked whether this was always the case. A counterexample is announced in \cite{BeKi} providing the existence of automorphisms of rational surfaces that do not leave invariant any curve.

\end{document}